\documentclass[12pt,a4,rm]{article}
\usepackage{amsfonts,amsthm,amssymb,amsmath,color,epsfig}
\usepackage{graphicx} 
\usepackage{xcolor}
\usepackage{setspace}
\definecolor{shinyblue}{rgb}{0, 0.3175, 0.8413}
\definecolor{shinygreen}{rgb}{0, 0.9527, 0.5238}
\definecolor{matlabyellow}{rgb}{0.9290, 0.6940, 0.1250}
\newcommand{\R}{{\mathbb R}}
\newcommand{\C}{{\mathbb C}}

\newcommand{\calP}{{\mathcal P}}
\newcommand{\Omegabar}{{\overline{\Omega}}}
\newcommand{\omegabar}{{\bar{\omega}}}
\newcommand{\pl}{{\!+\! }}
\newcommand{\mi}{{\!-\! }}
\newcommand{\rmv}{{\rm v}}
\newcommand{\ctilde}{{\tilde{c}}}
\newcommand{\gammabar}{{\bar{\gamma}}}
\newtheorem{theorem}{Theorem}[section]
\newtheorem{proposition}[theorem]{Proposition}

\newtheorem{lemma}[theorem]{Lemma}

\newtheorem*{theorem*}{Theorem}
\title{\Large Mass splitting in the time-discrete generalized Euler equations 
\\ and non-Monge solutions in multi-marginal optimal transport
}
\author{Gero Friesecke}
\date{January 5, 2026}

\begin{document}

\maketitle

\begin{small}
{\bf Abstract.} The time-discretized, spatially continuous generalized Euler equations are a prototype example of multi-marginal optimal transport, yet the question whether they exhibit mass-splitting (or equivalently, whether they have solutions that are not of Monge form) has remained open. Here we resolve this question by giving a mass-splitting example in one spatial dimension. Moreover we present a related and very simple fully discrete example of mass-splitting which reveals a transparent underlying mechanism. 

\end{small}
\vspace*{2mm}

\section{Introduction}

Arnold \cite{Ar66} made the celebrated observation that solutions to the incompressible Euler equations of fluid dynamics correspond to geodesics in the group of volume-preserving diffeomorphisms. The corresponding variational principle turns out to be ill-posed in general, which led Brenier \cite{Bre89} to introduce a relaxation which he showed to be well-posed. Physically this formulation, known as the \textit{generalized Euler equations}, allows mass splitting: a fluid particle can move from point A to point B via an ensemble of trajectories. Such mass-splitting phenomena are well-known in optimal transport. In fact, after time-discretization, 
\begin{itemize}
    \item Brenier's relaxation of Euler is a multi-marginal optimal transport (MMOT) problem in the Kantorovich formulation, where the incompressibility of the fluid is re\-presented by a marginal condition at each timepoint
    \item Arnold's variational principle is the Monge formulation of this MMOT problem.
\end{itemize}
See my recent textbook \cite{Fri25} for a comprehensive introduction to optimal transport  (including MMOT and important examples like the Euler equations).

In this paper, after briefly reviewing the different formulations of the Euler equations 
and their 
optimal transport interpretation, we prove that \\
- mass splitting still occurs in Brenier's relaxation after time-discretization, i.e. there \\
\textcolor{white}{-} exist non-Monge solutions \\
- mass splitting still occurs after also discretizing space. \\
In the light of previous such examples in continuous time and space \cite{Bre89, BFS09} (corresponding to the limit of infinitely many marginals) these results might not come as a surprise. Nevertheless, despite a considerable body of literature on mass-splitting in other MMOT problems (reviewed in section \ref{sec:review}) the case of time-discrete generalized Euler had remained open; the additional difficulty is that the marginals cannot be conveniently chosen, but are fixed to be the uniform measure. Moreover our fully discrete example is very  simple and reveals a transparent \textit{mechanism} leading to mass-splitting. 

The analytical examples reported here are in part  motivated by joint work in progress with Maximilian Penka on the numerical computation of solutions to the time-discrete generalized Euler equations \cite{FP26} via the algorithm introduced in \cite{FP23}. 


We close with a brief discussion from a modeling point of view: what is physically more correct, Euler (no mass splitting) or Brenier (mass splitting)?

It is a pleasure to dedicate this article to Willi J\"ager on the occasion of his 85$^{th}$ birthday. His deep understanding and 
practicing of applied mathematics all 
across 
inter\-disciplinary modeling, rigorous analysis, and numerical simulation 
continues to inspire.
    
\section{Different formulations of the Euler equations}
\subsection{Euler's formulation} 
Euler could build upon the work of distinguished researchers studying fluids before him: Archime\-des, Torricelli, Daniel Bernoulli, ... But he was the first to propose, in 1757, a complete model, by a system of partial differential equations:
\begin{align}
&   \partial_t ~\!u + (u\cdot \nabla)u=-\nabla p \;\;\; {\rm in } \;\Omega\times[0,T] \label{Euler} \\
&   {\rm div}\, u = 0 \; {\rm in} \; \Omega\times[0,T], \;\;\;\;\;\;\;\;\;\; u\cdot n = 0 \; {\rm on}\; \partial \Omega\times[0,T]. \label{Euler2}
\end{align}
Physically, $\Omega$ is the spatial region occupied by the fluid, $u$ is the velocity field, and $p$ is the pressure. 
Mathematically, $\Omega$ is an open bounded sufficiently regular (say, Lipschitz) domain in $\R^d$, 
$n$ is the outward unit normal to $\partial \Omega$, $u \, : \, \Omega\times[0,T]\to \R^d$ is a time-dependent vector field, $p\, : \, \Omega\times[0,T]\to{\mathbb R}$ is a time-dependent scalar field, and $u\cdot \nabla$ has the customary meaning $\sum_{i=1}^d u_i\partial_{x_i}$. The condition ${\rm div}\, u=0$ has the important meaning that the associated flow (see eq.~\eqref{flowmap}) is volume-preserving, that is to say the fluid is incompressible. Typically, one also prescribes the initial velocity, 
$$
   u|_{t=0} = u_0
$$
for some $u_0\, : \, \Omega\to\R^d$. 

The Euler equation \eqref{Euler} is the vanishing viscosity limit $\nu\to 0$ of the Navier-Stokes equation
\begin{equation*} 
    \partial_t ~\!u + (u\cdot \nabla)u= \nu \Delta u -\nabla p \;\;\; {\rm in } \;\Omega\times[0,T];
\end{equation*}
despite its neglect of viscous contributions to the stress it retains an important role in the description of near-inviscid phenomena such as turbulence (see e.g. \cite{DS13}).  

\subsection{Arnold's formulation}
Let us pass to a Lagrangian viewpoint, i.e. describe the flow by the position $g(x,t)$ at time $t$ of the fluid particle initially at $x$. The flow map $g(\cdot, t) \, : \, \Omega\to\Omega$ is defined by the ordinary differential equation
\begin{equation} \label{flowmap}
\begin{cases}
   \frac{d}{dt} g(x,t) = u( g(x,t),\,t) & \\
   \;\;\,\, g(x,0) = id. &
\end{cases}
\end{equation}
Since $u(\cdot,t)$ is, by eq.~\eqref{Euler2}, divergence-free, $g(\cdot,t)$ is volume-preserving, that is to say
\begin{equation} \label{det}
     \det Dg(x,t) \equiv 1 \; \forall x\in \Omega, \; \forall t\in[0,T]
\end{equation}
or equivalently  
$$
    g(\cdot,t)_\sharp 1_\Omega = 1_\Omega \; \forall t\in [0,T],
$$
where $T_\sharp \mu$ denotes the push-forward of a measure $\mu$ under a map $T$ (defined as $T_\sharp\mu(A)=\mu(T^{-1}(A))$ for any measurable set $A$) and $1_\Omega$ denotes the uniform measure on $\Omega$. This follows, e.g., from the elementary identity 
$$
 \frac{d}{dt} \det Dg(x,t) = {\rm div}\, u (y,t)\big|_{y=g(x,t)} {\rm det}\, Dg(x,t)
$$
for the flow map of any time-dependent vector field $u$. 

Arnold noticed that in Lagrangian coordinates, the Euler equations are (formally) equivalent to the variational principle
\begin{align} \label{arnold}
   {\rm min \; } \int_0^T  \big|\!\,\big|\tfrac{d}{dt}g(\cdot,t)\big|\!\,\big|_{L^2(\Omega)}^2 \; dt \;\; {\rm subject \; to} \; \, g(\cdot,t)_\sharp1_\Omega = 1_\Omega \; \forall t \in[0,T]
\end{align}
together with endpoint conditions
\begin{equation} \label{endpoint}
   g(\cdot,t)|_{t=0} = id, \;\;\; g(\cdot,t)|_{t=T} = g_*,
\end{equation}
where $g_* \, : \, \Omega\to\Omega$ is some prescribed volume-preserving map. The endpoint condition at time $T$ replaces the initial condition $u|_{t=0}=u_0$ for Euler. Informally, 
the variational principle says that 
$$
  \mbox{\it solutions to Euler = }L^2\, \mbox{\it geodesics in the group of volume-preserving maps.}
$$
For self-containedness we include a brief formal derivation of why (sufficiently smooth) minimizers of \eqref{arnold} are solutions to Euler. Introduce the set of volume-preserving maps,
\begin{align*}
    S = \{ \tilde{g} : \Omega \to \Omega \, | \, \tilde{g} \,{\rm measurable,} \; \tilde{g}_\sharp 1_\Omega \!=\! 1_\Omega \}.
\end{align*}
The abstract Euler-Lagrange equation of the variational principle is
\begin{align} \label{abstractEL}
   \tfrac{d^2}{dt^2} g(\cdot,t) \in \bigl(T_{g(\cdot,t)}S\bigr)^{\perp_{_{L^2}}}
\end{align}
where $T_{g(\cdot,t)}S$ denotes the tangent space of $S$ at the point $g(\cdot, )$ and $\perp_{L^2}$ denotes the orthogonal complement with respect to the $L^2$ inner product. It is not difficult to see that the tangent space is
given explicitly by
\begin{align*}
  T_{g(\cdot,t)}S = \{ {\rm v}(g(\!\,\cdot\!\,,t)) : \Omega\to\Omega \; | \,  \, {\rm div}\, {\rm v}\! =\! 0, \; {\rm v}\cdot \nu|_{\partial\Omega} \!=\! 0\}.
\end{align*}
We claim that the orthogonal complement with respect to the $L^2$ inner product is
\begin{align*}
   \bigl(T_{g(\cdot,t)}S\Bigr)^{\perp_{L^2}}  = \{\nabla p (g(\cdot,t)) \, | \, p :\Omega\to{\mathbb R}  \}.
\end{align*}
This follows from the calculation
\begin{align*}
   \bigl\langle \nabla p(g(\!\,\cdot\!\, ,t)),{\rm v}(g(\!\,\cdot\!\, ,t))\bigr\rangle_{L^2}
 & \; = \; \int_\Omega \nabla p(g(x,t)) \cdot {\rm v}(g(x,t)) \, dx  \\ & \!\!\! \underset{g(\cdot,t)\in S}{=} \int_\Omega \nabla p(y) \cdot {\rm n}(y) \, dy = \int_{\partial\Omega} p{\rm v}\cdot {\rm n} \, dS - \int_\Omega p \, {\rm div}\, {\rm v}
\end{align*}
where we have used the change of variables $y=g(x,t)$ and the fact that $g(\cdot,t)$ is volume-preserving (eq.~\eqref{det}). With this description of the orthogonal complement, the abstract Euler-Lagrange equation \eqref{abstractEL} becomes 
\begin{align} \label{concreteEL}
  \tfrac{d^2}{dt^2} g(\cdot,t) = -\nabla p(g(\cdot,t))\; {\rm for \; some \;} p \, : \, \Omega\times[0,T]\to\R.
\end{align}
Finally, let us compute the left hand side of \eqref{concreteEL}: by the chain rule, 
\begin{align*}
   \tfrac{d^2}{dt^2} g(\!\,\cdot\!\, ,t) = \tfrac{d}{dt} u(g(\!\, \cdot \!\, ,t),t) = \bigl(\partial_t u + \underbrace{Du \, u}_{=\, (u\cdot\nabla) u} \!\bigr)(g(\!\,\cdot\!\, ,t),t),
   \vspace*{-2mm}
\end{align*}
so \eqref{concreteEL} is precisely the Euler equation \eqref{Euler}. 
\subsection{Brenier's formulation} \label{sec:GenEuler}
Ebin and Marsden \cite{EM70} showed that 
Arnold's variational principle has a unique optimizer when the endpoint map $g_*$ is a smooth volume-preserving diffeomorphism which is sufficiently close to the identity in a suitable Sobolev norm. By contrast, Shnirelman \cite{Sh87} proved the somewhat surprising result that there exist smooth volume-presering diffeomorphisms $g_*$ on $\Omega=[0,1]^3$ for which the variational principle has no minimizer.    

This ill-posedness led Brenier \cite{Bre89} to introduce a very interesting relaxation. First one re-writes Arnold's variational principle via Fubini's theorem as 
\begin{align}
  & {\rm min \; } \int_\Omega \int_0^T \big|\tfrac{d}{dt}g(x,t)\big|^2 dt \, dx \;\; {\rm subject \; to} \; \, g(\cdot,t)_\sharp1_\Omega = 1_\Omega \; \forall t \\
  & {\rm with \; endpoint \; condition \;}g(\cdot,0)=id, \;\; g(\cdot, T)=g_*.
\end{align}
The inner integral is the action of the path taken by the fluid particle initially at $x$. 

\begin{figure}[http!]
    \begin{center}
        \includegraphics[width=0.75\textwidth]{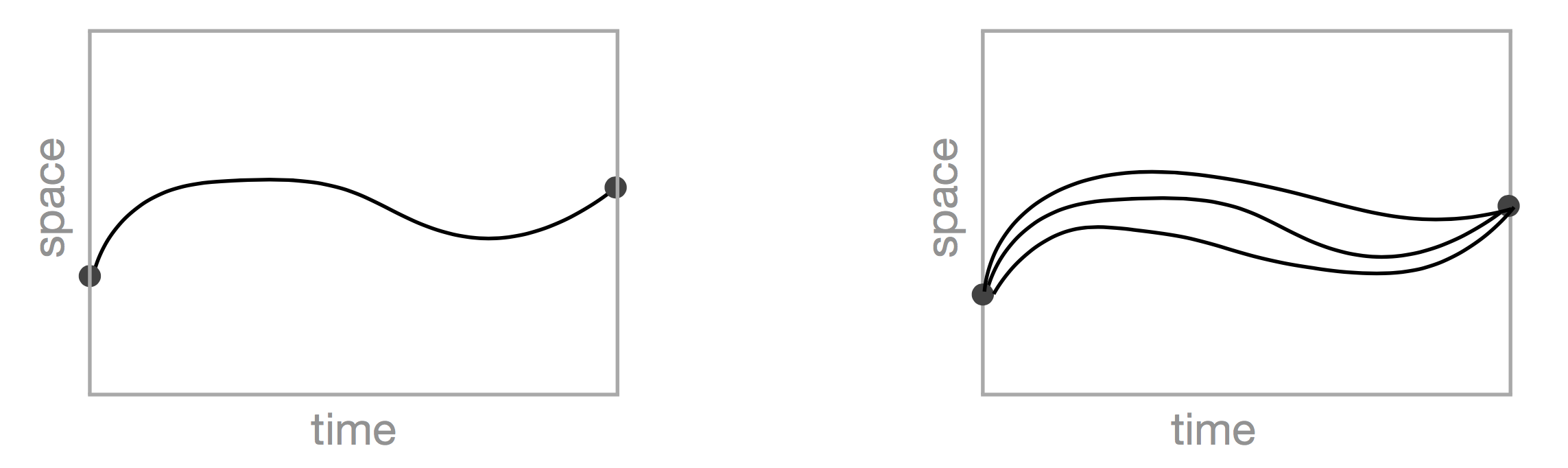}
        \vspace*{-4mm}
        
    \end{center}
    \caption{\textit{Left}: each fluid particle follows a single path $\omega$ (Arnold's variational principle). \textit{Right}: each fluid particle can follow multiple paths (Brenier's variational principle).}
    \label{F:paths}
\end{figure}

Now instead of each fluid particle following a single path, one allows it to follow an ensemble of paths described by a probability measure on path space, 
\begin{align} \label{brenier1}
&   {\rm min}\, \int_{\omega \in C([0,T];\Omega)} \int_0^T |\dot{\omega}(t)|^2 dt \, d\gamma(\omega)  \;\;\; {\rm subject \; to} \; \pi_t{}_\sharp \gamma = {\mathbf 1}_\Omega\;\forall t \\ 
& {\rm with \; endpoint \; condition \;}
(\pi_0,\pi_T)_\sharp\gamma = (id,g_*)_\sharp {\mathbf 1}_\Omega.
\label{brenier2}
\end{align}
Here the minimization is over probability measures on path space, $\gamma \in {\mathcal P}(C([0,T];\Omega))$, and 
\vspace*{-4mm}
\begin{align*}
     \pi_t \, : \, C([0,T];\Omega)& \to\Omega \\
        \omega & \mapsto \omega(t)
\end{align*}
denotes the map which assigns to each path its value at time $t$. As before, $g_*\, : \, \Omega\to\Omega$ is a given volume-preserving map. The variational problem 
\eqref{brenier1} is known as \textit{generalized Euler}. As shown by Brenier, with respect to the narrow topology on path space the new functional is lower semi-continuous and its domain is compact, ensuring existence of optimizers. 
Physically, the passage from single paths for each fluid particle to ensembles of paths means that one allows mass splitting. See Figure \ref{F:paths}. 

There are endpoint conditions which lead to mass splitting. 
\\[2mm]
{\bf Example 1.} This example in one dimenson is due to Brenier \cite{Bre89}: $\Omega=[-1,1]$, $g_*(x)=-x$ (i.e. one turns the fluid upside down), $T=\pi$. In this case Arnold's variational principle has no solution, and Brenier's variational principle is uniquely solved by a certain probability measure $\gamma$ which is concentrated for each initial position $x\in[-1,1]$ on the one-parameter family of paths
$$
   \omega_{x,\rmv}(t) = g_{\rmv}(x,t) = x \cos t + 
\rmv \sin t  \;\; \;\;(\rmv\in [-\sqrt{1-x^2}, \sqrt{1-x^2}])
$$
all of which connect the initial point $x$ to the endpoint $-x$. This form of paths arises by solving  eq.~\eqref{concreteEL} with the time-independent pressure $p(x)=\tfrac12 x^2$. More precisely, $\gamma$ is concentrated on these paths with probability density
$$
    f(x,\rmv) = \frac{1}{2\pi\sqrt{1 - (x^2+\rmv^2)}}.
$$
For a visualization of this example see \cite{Fri25} Figure 1.19. 
\\[2mm]
{\bf Example 2.} Building upon Example 1, Bernot, Figalli and Santambrogio \cite{BFS09} constructed (not necessarily unique) mass-splitting solutions in the two-dimensional disc subject to the endpoint condition $g_*$ being a rotation. 
\\[2mm]
One might argue that mass splitting indicates some sort of breakdown of the original Euler model of fluid dynamics. For further discussion of this point from a modeling point of view see section \ref{sec:modeling}. 

\subsection{Time-discretization and interpretation as optimal transport} \label{sec:TD}

Replacing the continuous time interval $[0,T]$ by a discrete set of times $\{t_0,t_1,...,t_N\}$ with $0=t_0<t_1<...<t_N=T$ leads to the following simplifications: 
\begin{eqnarray*}
    \mbox{path space:} \; C([0,T];\Omega)  & \;\; \rightsquigarrow \;\; & \Omega^{N+1} \\
    \mbox{paths:} \; \omega\in C([0,T];\Omega) & \;\; \rightsquigarrow \;\; & (\omega_0,...,\omega_N)\in \Omega^{N+1} \\
    \mbox{probability measures: } \gamma\in \calP(C([0,T];\Omega)) & 
    \;\; \rightsquigarrow \;\; &
    \gamma\in \calP(\Omega^{N+1}).
\end{eqnarray*}
Brenier's variational principle \eqref{brenier1} reduces to: 
\begin{align} \label{MMOT'1}
&   \min_{\gamma\in\calP(\Omega^{N+1})}\, \int_{\Omega^{N+1}} \sum_{i=1}^N \frac{|\omega_i-\omega_{i-1}|^2}{t_i-t_{i-1}} \, d\gamma(\omega_0,...,\omega_N) \;\; {\rm subject \; to} \; \pi_i{}_\sharp \gamma = {\mathbf 1}_\Omega\,\forall i=0,...,N \\
& \;\,{\rm with \; endpoint \; condition}\; (\pi_0,\pi_N)_\sharp \gamma = (id,g_*)_\sharp {\mathbf 1}_\Omega. \label{MMOT'2}
\end{align}
This discrete version was also introduced by Brenier \cite{Bre93}. 
Conversely, as shown by Nenna \cite{Nen17} building on earlier results in \cite{Bre93}, the discrete problem Gamma-converges in a suitable sense to the continuous problem as the time stepsize goes to zero. 

A further simplification can be achieved by subsuming the endpoint condition \eqref{MMOT'2} into the cost. This condition implies that any path  $(\omega_0,...,\omega_N)$ charged by an admissible competitor $\gamma$ must satisfy $\omega_N=g_*(\omega_0)$, and is thus completely characterized by its first $N\!-\!1$ components $(\omega_0,...,\omega_{N-1})$, leading to the following reduced variational principle (see Lemma \ref{L:equiv} below): 
\begin{align}
&   \min_{\gammabar\in\calP(\Omega^{N})}\, \int_{\Omega^{N}} \Bigl(\sum_{i=1}^{N-1} \frac{|\omega_i-\omega_{i-1}|^2}{t_i-t_{i-1}} + \frac{|\omega_{N-1}-g_*(\omega_0)|^2}{t_N-t_{N-1}}\Bigr)\, d\gammabar(\omega_0,...,\omega_{N-1}) \label{MMOT1} \\
& {\rm subject \; to} \; \pi_i{}_\sharp \gammabar = {\mathbf 1}_\Omega\,\forall i=0,...,N\,\mi\, 1. \label{MMOT2}
\end{align}
Problems \eqref{MMOT'1}--\eqref{MMOT'2} and \eqref{MMOT1}--\eqref{MMOT2} are instances of multi-marginal optimal transport (MMOT) problems 
\begin{align} \label{Kant1}
&   \min_{\gamma\in\calP(X_1\times ... \times X_N)}\, \int_{X_1\times ...\times X_N} c(x_1,...,x_N)\, d\gamma(x_1,...,x_N) \\
& {\rm subject \; to } \; \pi_i{}_\sharp \gamma = \mu_i \label{Kant2}
\end{align}
where 
\begin{itemize}
\vspace*{-2mm}

    \item the $X_i$ are metric spaces
    \vspace*{-2mm}
    
    \item $\pi_i$ is the projection map of the product space $X_1\times ... \times X_N$ onto the $i$-th factor (i.e. $\pi_i(x_1,...,x_N)=x_i$) and the $\mu_i$ are prescribed Borel probability measures on the factor spaces $X_i$, and so the constraint \eqref{Kant2} amounts to prescribing the $N$ marginals of the multivariate probability measure $\gamma$
    \vspace*{-2mm}
    
    \item $c\, : \, X_1\times ...\times X_N\to\R$ is a cost function. 
\end{itemize}
\vspace*{-2mm}

\noindent
Thus one seeks to minimize the expected value of some cost function $c$ over multivariate probability measures $\gamma$ with prescribed marginals. Minimizers are called optimal plans, and can be proven to exist in great generality (for instance when the $X_i$ are arbitrary closed subsets of $\R^d$, the $\mu_i$ are arbitrary Borel probability measures on $X_i$, and $c$ is lower semi-continuous and bounded from below (\cite{Fri25} Theorem 3.1)). 
Such problems arise in numerous other contexts such as interpolation of data in time, interpolation of data in space, or many-electron quantum mechanics, with different applications corresponding to different cost functions $c$; see \cite{Fri25} Section 1.6. 

For future use let us make the reduction to  \eqref{MMOT1}--\eqref{MMOT2} precise.

\begin{lemma} \label{L:equiv} For any closed and bounded  subset $\Omega\subset\R^d$ and any $N\ge 2$, a probability measure $\gamma\in\calP(\Omega^{N+1})$ solves \eqref{MMOT'1}--\eqref{MMOT'2} if and only if its marginal $(\pi_0,...,\pi_{N-1})_\sharp \gamma\in\calP(\Omega^N)$ with respect to the first $N$ copies of $\Omega$ solves \eqref{MMOT1}--\eqref{MMOT2}. 
\end{lemma}

{\bf Proof} Let $\gamma$ be an admissible competitor of \eqref{MMOT'1}--\eqref{MMOT'2}. Its marginal  $\gammabar =(\pi_0,...,\pi_{N-1})_\sharp \gamma$ with respect to the first $N$ copies of $\Omega$ is then admissible for \eqref{MMOT1}--\eqref{MMOT2}. Moreover by \eqref{MMOT'2}, any path $(\omega_0,...,\omega_N)$ in the support of $\gamma$ must satisfy $\omega_N=g_*(\omega_0)$ and hence
$$
\int_{\Omega^{N+1}} \underbrace{\sum_{i=1}^N \frac{|\omega_i-\omega_{i-1}|^2}{t_i-t_{i-1}}}_{=:c(\omega_0,...,\omega_N)} d\gamma = \int_{\Omega^{N+1}}  \Bigl(\underbrace{\sum_{i=1}^{N-1} \frac{|\omega_i-\omega_{i-1}|^2}{t_i-t_{i-1}} + \frac{|\omega_{N-1}-g_*(\omega_0)|^2}{t_N-t_{N-1}}}_{=:c_*(\omega_0,...,\omega_{N-1})}\Bigr)\, d\gamma = \int_{\Omega^N} c_* d\gammabar. 
$$
Conversely, any admissible competitor $\gammabar\in\calP(\Omega^N)$ of \eqref{MMOT1}--\eqref{MMOT2} can be extended to an admissible competitor of \eqref{MMOT'1}--\eqref{MMOT'2}, namely $\gamma=e_\sharp \gammabar$ where $e$ is the extension map $e(\omega_0,...,\omega_{N-1})=(\omega_0,...,\omega_{N-1},g_*(\omega_0))$, and we have by the change-of-variables formula 
$$
   \int_{\Omega^N} c_* \, d\gammabar = \int_{\Omega^N} c\circ e \; d\gammabar = \int_{\Omega^{N+1}} c \; d e_{\sharp}\gammabar = \int_{\Omega^{N+1}} c \; d\gamma.
$$

\section{Mass-splitting in optimal transport}
\label{sec:splitting} An important, longstanding open problem in multi-marginal optimal transport is to understand when mass splitting occurs, or equivalently when optimal plans are not of ``Monge'' form. Let us formalize these notions
and summarize previous results.
\subsection{Terminology; Monge versus Kantorovich} 
{\bf Definition} {\it An optimal plan for the general $N$-marginal optimal transport problem \eqref{Kant1}--\eqref{Kant2} is called \emph{mass-splitting with respect to the $i$-th marginal} if 
$\gamma$ gives mass to different configurations $(x_1,...,x_N)$ with same $x_i$, and \emph{everywhere mass-splitting} if it is mass-splitting with respect to every marginal.
   
In particular, an optimal plan for the time-discretized generalized Euler equations \eqref{MMOT1}--\eqref{MMOT2} is called mass-splitting at the discrete time $t_i$ ($i\in\{0,...,N\!-\! 1\}$) if $\gamma$ gives mass to different paths $(\omega_0,...,\omega_N)$ with same $\omega_i$, and everywhere mass splitting if it is mass-splitting for every discrete time $t_i$.  
} 
\\[2mm]
Thus an optimal plan is mass-splitting with respect to the $i$-th marginal (or at the discrete time $t_i$) if and only if it is not of Monge form with respect to the $i$-th marginal (or at the discrete time $t_i$), where we recall the following standard definition: 
\\[2mm]
{\bf Definition} (see \cite{Fri25} Section 1.6) {\it An optimal plan for the general $N$-marginal optimal transport problem \eqref{Kant1}--\eqref{Kant2} is \emph{of Monge form with respect to the $i$-th marginal} if there exist measurable maps $T_k \, : \, X_i\to X_k$ ($k\in\{1,...,N\}\backslash\{i\}$) such that 
$$
  \gamma = (T_1,...,T_{i-1},id,T_{i+1},...,T_N)_\sharp \mu_i.
$$
} 

\noindent
It is clear from these definitions that any non-mass-splitting optimal plan is a solution to the Monge formulation of the multi-marginal optimal transport problem \eqref{Kant1}--\eqref{Kant2},
\begin{align}
     \label{Monge1}
&   \min_{T_1,..,T_{i-1},T_{i+1},...,T_N}\, \int_{X_i} c(T_1(x_i),...,T_{i-1}(x_i),x_i,T_{i+1}(x_i),...,T_N(x_i))\, d\mu_i(x_i) \\
& {\rm over \; measurable \; maps \; }T_k \, : \, X_i\to X_k {\rm \; subject \; to } \; T_k{}_\sharp \mu_i = \mu_k \; (k\in\{1,...,N\}\backslash\{i\}). \label{Monge2}
\end{align} 
%


\subsection{Previous results} \label{sec:review}
For $N=2$ and in Euclidean spaces (i.e. $X_1$, $X_2$ are closed subsets of $\R^d$) the question of mass-splitting is reasonably well understood (see e.g. \cite{Fri25}): if the cost function $c$ is differentiable and satisfies the twist condition
$$
   y \mapsto \nabla_x c(x,y) \; {\rm injective},
$$
and the first marginal $\mu_1$ is absolutely continuous with respect to the Lebesgue measure, then optimal plans are unique and of Monge form with respect to the first marginal, that is to say no mass-splitting occurs. The celebrated special case $c(x,y)=|x-y|^2$, which is due to Brenier \cite{Bre91}, was in fact motivated by his seeking to understand what happens for the generalized Euler equations in a single timestep, and can be formulated as follows:
\begin{theorem*}{\rm (Brenier's theorem \cite{Bre91})} 
Let $\Omega$ be an open bounded domain in $\R^d$, and let $\mu_1$, $\mu_2$ be absolutely continuous probability measures in $\calP(\Omegabar)$. For the problem 
$$
  \min_{\gamma\in\calP(\Omegabar\times\Omegabar)} \int _{\Omegabar\times\Omegabar}|\omega_0-\omega_1|^2 d\gamma(\omega_0,\omega_1) \; \; 
  {\rm subject\; to \; \;}\pi_i{}_\sharp \gamma=\mu_i \; (i=1,2)
$$
no mass-splitting occurs, that is to say the optimal $\gamma$ is unique and does not give mass to different paths $(\omega_0,\omega_1)$ with same $\omega_0$, nor different paths $(\omega_0,\omega_1)$ same $\omega_1$.  
\end{theorem*}

For $N\ge 3$ the situation is much less well understood, except for certain examples. 

Agueh and Carlier \cite{AC11} showed that there is no mass-splitting and optimal plans are unique for the Wasserstein barycenter cost 
\begin{align}
   c(x_1,...,x_N) = \sum_{i=1}^N \lambda_i|x_i - B(x)|^2, \; B(x)=\sum_{i=1}^N \lambda_i x_i
\end{align}
arising in the spatial interpolation of data. Here 
the $\lambda_i$ are positive weights with $\sum_{i=1}^N \lambda_i=1$. In the prototypical case of equal weights this follows from earlier results by Gangbo and \'Swiech \cite{GS98}. For extensions to mathematically similar costs arising in economics (multi-agent matching) and data interpolation ($p$-Wasserstein barycenters) see \cite{Pas14, BFR24}. 

Pass \cite{Pas13} exhibited an example of a (non-unique) mass-splitting optimal plan for the Coulomb cost 
\begin{align}
   c(x_1,...,x_N) = \sum_{1\le i<j\le N}\frac{1}{|x_{i}- x_j|}
\end{align}
arising in electronic structure, in dimension $d=3$ and for $N=3$ marginals. It is not known for this cost whether there always exists an optimal Monge plan, except when $d=1$ and $N$ is arbitrary \cite{CDD15}, in which case the answer is Yes. 

Friesecke \cite{Fri19} gave an example of a unique mass-splitting optimal plan for the Frenkel-Kontorova-type cost 
\begin{align}
   c(x_1,...,x_N) = \sum_{1\le i<j\le N} v(|x_{i}- x_j|), \;\; v(r)=\frac{r^4}{4}-\frac{r^3}{3}
\end{align}
arising in statistical mechanics, in dimension $d=1$ and for $N=3$ uniform marginals. For the analogous cost with the repulsive harmonic potential $v(r)=-r^2$, Gerolin, Kausamo and Rajala proved in general dimension $d$ and for $N=3$ that there are absolutely continuous marginals for which mass-splitting must occur (i.e. the Monge problem has no solution); examples of (not necessarily unique) mass-splitting plans for this cost had been found earlier by Pass \cite{Pas12}. 

Benamou, Gallouet, and Vialard \cite{BGV19} found an example of a (not necessarily unique) mass-splitting optimal plan for the Wasserstein cubic spline cost
\begin{align}
   c(x_1,...,x_N) = \sum_{i=2}^{N-1} |x_{i-1}-2x_i+x_{i+1}|^2
\end{align}
arising in the time-interplation of data,
in dimension $d=1$ and with $N=3$ marginals one of which is not absolutely continuous. The authors also proved that in their example mass-splitting must occur (i.e. the Monge problem has no minimizer). It is not known whether these findings persist if all marginals are required to be absolutely continuous. 

A nontrivial upper bound on the support dimension of optimal plans for general costs was obtained by Pass \cite{Pas12}, using techniques from geometric measure theory. 


\section{Why does the mass want to split in the generalized Euler equations? A simple mechanism}

Previous understanding \cite{Bre89, BFS09} requires a deep analysis of the full continuous path space model \eqref{brenier1}--\eqref{brenier2} and rests on indirect arguments: after ingeniously guessing mass-splitting minimizers one verifies their optimality with the help of Kantorovich duality.

Here we pass to the simplified setting of discrete space and time, and give \\
(1) a very simple discrete 
mass-splitting example (see Figure 2)
\\
(2) a transparent variational argument why it is favourable over non-splitting solutions \\
(3) a fairly short variational argument why it is optimal. \\
While coming up with the example was of course inspired by Brenier's pioneering work, the advance lies in the variational argument (2), which reveals a transparent \textit{mechanism} leading to mass-splitting. A verbal summary of the mechanism is given after the proof of Proposition \ref{P1}~a). 

We use the (minimalistic) space discretization
$$
   [-1,\, 1] \; \rightsquigarrow \; \Omega=\{-1,0,1\}
$$
and a uniform time discretization with any number of timesteps
$$
  [0,\,T] \; \rightsquigarrow \{t_0,...,t_N\}, \; t_j=j,
$$
as well as Brenier's endpoint map
\vspace*{-2mm}
$$
  g_*(\omega_0) = -\omega_0 
$$
(i.e. one turns the fluid upside down). 
Candidate or optimal plans $\gamma\in\calP(\Omega^{N+1})$ can then be indentified with their density with respect to counting measure, i.e. they can be viewed as functions $\gamma \, : \, \Omega^{N+1}\to\R$ on paths $\omega=(\omega_0,...,\omega_{N})$ satisfying  $\gamma\ge 0$ and $\sum_{\omega\in \Omega^{N+1}} \gamma(\omega) = 1$. Moreover, as explained in section \ref{sec:TD}, due to the endpoint condition $\omega_N=-\omega_0$ candidate or optimal plans are supported on the set of paths  $\{\omega\in\{-1,0,1\}^{N+1} \, : \, \omega_N=-\omega_0\}$ and the problem  \eqref{MMOT'1}--\eqref{MMOT'2} becomes
\vspace*{-4mm}
\begin{align} \label{discbren1}
    &   \min_{\gamma\in\calP\bigl(\{\omega\in\{-1,0,1\}^{N+1}\, : \, \omega_N=-\omega_0\}\bigr)}  \, \sum_{\substack{\omega\in\{-1,0,1\}^{N+1} : \\ \omega_N=-\omega_0}} \Bigl(\sum_{i=1}^{N} |\omega_{i-1}\!-\!\omega_{i}|^2  \Bigr) \; \gamma(\omega) \\
    & \;\;\, {\rm subject \; to \; } \; \pi_i{}_\sharp \gamma = \tfrac{1}{3}1_\Omega\; \forall i. \label{discbren2}
\end{align}
\begin{proposition} \label{P1} {\rm (Mass-splitting in the fully discrete case)} Let $N\ge 3$.\\[1mm]  
{\rm a)} For $N\ge 4$, any minimizer of \eqref{discbren1}--\eqref{discbren2} is mass-splitting. 
\\[1mm]
{\rm b)} A minimizer of \eqref{discbren1}--\eqref{discbren2} is given by the probability measure concentrated on the paths shown in Figure 2 with 
\vspace*{-2mm}
\begin{eqnarray*}
   & {\rm probability}\bigl( {\rm thick\; green\; path} \bigr) & \hspace*{-2.3mm} = \frac{N-3}{3(N-1)}, \\
  & {\rm probability}\bigl( {\rm any\; other\; path} \bigr) & \hspace*{-2.3mm} = \frac{1}{3(N-1)}.
\end{eqnarray*}
{\rm c)} This minimizer is unique if and only if $N$ is even. 
\end{proposition}
\begin{figure}[http!]
\begin{center}
\includegraphics[width=0.5\textwidth]{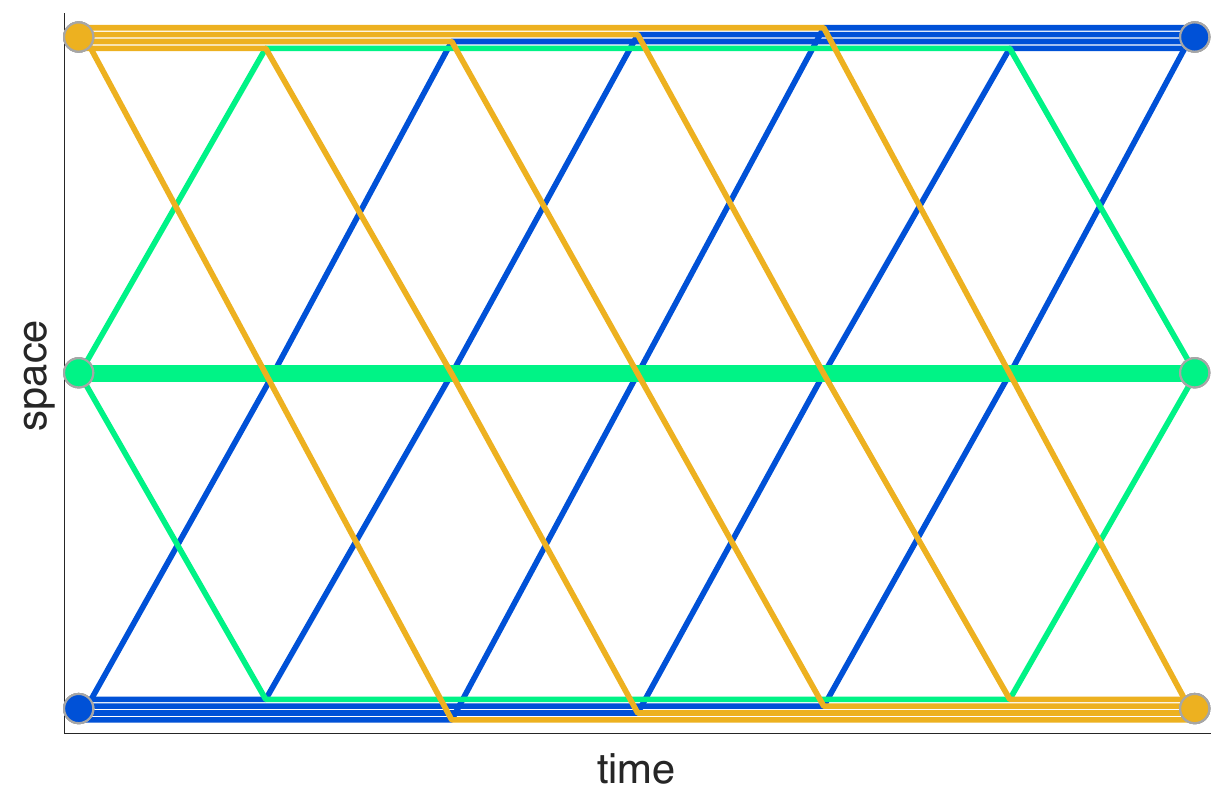}
\caption{Solution to the fully discrete generalized Euler equations \eqref{discbren1}--\eqref{discbren2}} 
\end{center}
\label{F:gerosexample}
\end{figure}

\noindent
The picture corresponds to the following formula for the optimizer $\gamma_0$, with the yellow, green, and blue paths 
representing the contributions in the first, second, and third line: 
\begin{align} 
 \gamma_0(\omega) &= \frac{1}{3(N-1)} \sum_{\nu=1}^{N-1} 
\Bigl(\prod_{i=0}^{\nu-1} \delta_1(\omega_i)\Bigr)
 \delta_0(\omega_\nu) 
 \Bigl(\prod_{j=\nu+1}^{N} \delta_{-1}(\omega_j)\Bigr) \nonumber \\
 &+ \frac{1}{3(N-1)} \Bigl[  (N-3) \Bigl(\prod_{i=0}^N \delta_0(\omega_i)\Bigr) + 
 \delta_0(\omega_0) \Bigl(\prod_{i=1}^{N-1} \delta_1(\omega_i) + \prod_{i=1}^{N-1}\delta_{-1}(\omega_i)\Bigr) \delta_0(\omega_N)  \Bigr] \nonumber \\
 &+ \frac{1}{3(N-1)} \sum_{\nu=1}^{N-1} 
 \Bigl(\prod_{i=0}^{\nu-1} \delta_{-1}(\omega_i)\Bigr)
 \delta_0(\omega_\nu) 
 \Bigl(\prod_{j=\nu+1}^{N} \delta_{1}(\omega_j)\Bigr). \label{gerosplan}
\end{align}

Before proving that the above plan is optimal, we give a simple explanation why the optimizer must necessarily exhibit mass-splitting.
\vspace*{2mm}

{\bf Proof of mass-splitting (Prop.~\ref{P1}~a)).} We denote the cost to be minimized by 
\begin{align*}
   C[\gamma] = \sum_{\substack{\omega\in\{-1,0,1\}^{N+1}: \\ \omega_N=-\omega_0}} \Bigl(\underbrace{|\omega_0-\omega_1|^2 + |\omega_1-\omega_2|^2 +...+ |\omega_{N-1}-\omega_N|^2}_{=:c(\omega)}\Bigr) \cdot \gamma(\omega).
\end{align*}

1. The outer orbits (i.e. those starting at $1$ or $-1$) must move.  Trivially, moving monotonically is cheaper than moving non-monotonically. The target position can be reached monotonically either by 2 small steps (i.e. $|\omega_{i}-\omega_{i-1}|=1$) or 1 large step (i.e. $|\omega_i-\omega_{i-1}|=2$). By strict convexity of the cost function $h(\omega_i-\omega_{i-1})=|\omega_i-\omega_{i-1}|^2$, 2 small steps are better than 1 large step: 
\begin{align}
  1^2  +  1^2   &< \; 2^2 \;\; + \;\; 0^2, \label{keyineq} \\
  {\rm or \; more \; generally}\;\; h(\tfrac{2+0}{2}) &< \frac{h(2) + h(0)}{2} \;\;  {\rm for \; any \; strictly \; convex \; }h. \nonumber
\end{align}
In particular, the minimum cost $c(\omega)$ for any outer orbit equals the left hand side of \eqref{keyineq}. 

2. If the middle orbit (i.e. the one starting at $0$) doesn’t move, the outer orbits must make large steps. This is because, by incompressibility (eq.~\eqref{discbren2}), the middle site $0$ is then fully occupied by the middle orbit at all times. This yields the following lower bound for the total cost:
\begin{align}
   C[\gamma]\; = \!\!\sum_{\substack{\omega \, : \, \omega_0=1, \\ \omega_N=-1}}\!\! c(\omega)\cdot\gamma(\omega) + \!\!\!\sum_{\substack{\omega \, : \, \omega_0=-1, \\ \omega_N=1}}\!\! c(\omega)\cdot\gamma(\omega) + \!\!\!\sum_{\substack{\omega\, : \, \omega_0= 0, \\ \omega_N=0}}\!\! c(\omega)\cdot\gamma(\omega) \ge  (\, 2^2 \;\; + \;\; 2^2 \;\; + \;\; 0^2 \,) \cdot \tfrac{1}{3}. \label{doesntmove}
\end{align}
   
3. If the middle orbit moves (and doesn't split), the outer orbits can make small steps, but the cost $c(\omega)$ for the middle orbit goes up from $0$ to  $1^2+1^2$.  (More precisely, when $N\ge 3$ and the middle orbit vacates the middle site for two times, the outer orbits can both make small steps.) This yields the lower bound
\begin{align}
   C[\gamma]\;  \ge  (2\!\cdot\! 1^2 + 2\!\cdot\!1^2 + 2\!\cdot\!1^2) \cdot \tfrac{1}{3} \label{moves}
\end{align}
which is achieved, e.g., when $\gamma$ is concentrated with equal probability on the three paths $(1,0,-1,...,-1)$, $(-1,-1,0,1,...,1)$, and $(0,1,1,0,...,0)$). Comparing with \eqref{doesntmove} shows that this is favourable over the middle orbit not moving. In total we conclude that the right hand side of \eqref{moves} is the optimal total cost when there is no splitting (i.e. when each set $\{\omega \, : \, \omega_0=x\}$ charged by $\gamma$ is a singleton). 

4. If the middle orbit moves only \textcolor{shinygreen}{partially} (i.e. splits into a part that moves and another that doesn't), as in Figure 2, the outer orbits can still make small steps, by sending one $(N\!-\!1)^{th}$ of their mass through the middle at each of the times $t_1,...,t_{N-1}$. This only requires the middle orbit to vacate two $(N\!-\!1)^{ths}$ of its mass at these times, so the middle orbit pays only $\tfrac{2}{N-1}(1^2+1^2)$, giving the total cost
\begin{align} \label{opticost}
   C[\gamma_0]\;  =  (2\!\cdot\! 1^2 + 2\!\cdot\!1^2 + \textcolor{shinygreen}{\tfrac{2}{N-1}} 2\!\cdot\!1^2) \cdot \tfrac{1}{3}.
\end{align}
This is lower than the optimal cost when there is no splitting (r.h.s. of \eqref{moves}) when $N\ge 4$.
\\[3mm]
Let us summarize the mechanism leading to mass-splitting as revealed by the proof:
\begin{itemize}
\vspace*{-2mm}
    \item \textit{By convexity of the Lagrangian in the velocity, outer orbits want to cover the required distance via small steps not large steps, and hence must pass through the midpoint.}
\vspace*{-2mm}
    \item \textit{By mass conservation, this is only possible if the middle orbit vacates the midpoint, at some cost.} 
\vspace*{-2mm}
    \item \textit{This cost can be greatly reduced when the middle orbit vacates the midpoint only partially and the outer orbits send their mass through the midpoint in small pieces at many different times.} 
\end{itemize}

\vspace*{2mm}

{\bf Proof of formula for optimizer (Prop.~\ref{P1}~b)).}
Given any candidate plan $\gamma$, it is useful to investigate its behavior on the following subsets of path space:
\begin{small}
\begin{align*}
\mbox{inner orbits} &:= \{\omega\in\{-1,0,1\}^{N+1} \, : \, \omega_0=0,\,\omega_N=0\}, \\
\mbox{outer orbits} &:= \{\omega\in\{-1,0,1\}^{N+1} \, : \, \omega_0=1,\,\omega_N=\!-1  \mbox{ or }\omega_0=1,\omega_N=\!-1\}, \\
\mbox{orbits that vacate 0} &:= \mbox{inner orbits s.th.~there exists }i\in\{1,...,N\mi 1\}\mbox{ with }\omega_i\neq 0, \\
\mbox{orbits that stays 0} &:= \mbox{inner orbit s.th.~}\omega_i=0\mbox{ for all }i, \\
\mbox{orbits that pass through 0} &:= \mbox{outer orbits s.th.~there exists }i\in\{1,...,N\mi 1\}\mbox{ with }\omega_i=0, \\
\mbox{orbits that don't pass through 0} &:= \mbox{outer orbits s.th.~}\omega_i\neq 0 \mbox{ for all }i.
\end{align*}
\end{small}

\noindent
Note that $\gamma$ is supported on the disjoint union of the last four sets, which equals the disjoint union of the first two sets.

The lower bounds derived in the proof of a) give 
\begin{align}
C[\gamma|_{\mbox{\rm \scriptsize inner \, orbits}}] & \ge (1^2 \pl 1^2) \cdot \gamma(\mbox{\scriptsize orbits that vacate 0}) \; + \; 0 \cdot \gamma(\mbox{\scriptsize orbits that stay 0}), \label{clb} \\
C[\gamma|_{\mbox{\rm \scriptsize outer\, orbits}}] & \ge 2^2 \cdot \gamma(\mbox{\scriptsize orbits that don't pass through 0}) \; + \; (1^2 \pl 1^2)\cdot \gamma(\mbox{\scriptsize orbits that pass through 0}). \nonumber
\end{align}
By mass conservation, eq.~\eqref{discbren2}, we have $\pi_i{}_\sharp\gamma (0) = \tfrac13$ for all $i$, so the mass of orbits that pass through $0$ at time $t_i$  equals that of orbits that vacate $0$ at time $t_i$. 
Since each orbit that passes through $0$ 
can pass through it 
at a single time $t_i$ but each orbit that  vacates $0$ can vacate it for a maximum of $N\!-\!1$ timepoints (namely at $t_1,...,t_{N-1}$), we have 
\begin{equation} \label{massineq}
     \gamma (\mbox{\footnotesize orbits that pass through 0}) \le (N-1) \cdot \gamma(\mbox{\footnotesize orbits that vacate 0}).
\end{equation}
Moreover since $\pi_0{}_{\sharp}\gamma=\tfrac13 1_{\Omega}$ we have 
\begin{equation} \label{masseq}
     \tfrac23 = \gamma\bigl(\mbox{\footnotesize outer  orbits}\bigr) = \gamma(\mbox{\footnotesize orbits that pass through 0}) + \gamma(\mbox{\footnotesize orbits that don't pass through 0}).
\end{equation}
Using first the lower cost bounds \eqref{clb} and then the mass inequality and equality  \eqref{massineq} and \eqref{masseq} gives
\begin{align}
    C[\gamma] &= C[\gamma|_{\mbox{\scriptsize inner  orbits}}] + C[\gamma|_{\mbox{\scriptsize outer \, orbits}}]\nonumber \\
    &\ge 2 \gamma(\mbox{\footnotesize orbits that vacate 0}) + 2  \gamma(\mbox{\footnotesize orbits that pass through 0})\nonumber \\
    & \;\;\;\; + 4 \gamma(\mbox{\footnotesize orbits that don't pass through 0})\nonumber \\
    &\ge \tfrac{2}{N-1} \gamma(\mbox{\footnotesize orbits that pass through 0}) + 2 \gamma(\mbox{\footnotesize orbits that pass through 0}) \nonumber \\
    & \;\;\;\; + 4 \bigl[\tfrac23 - \gamma(\mbox{\footnotesize orbits that pass through 0})\bigr] \nonumber \\
    &= \tfrac{2}{N-1}\alpha + 2\alpha + 4 (\tfrac23 - \alpha) \; = \; \tfrac83 + \bigl( \tfrac2{N-1} - 2\bigr)\alpha 
    \label{totcostbd}
\end{align}
where
$$
   \alpha := \gamma(\mbox{\footnotesize orbits that pass through 0}) \in [0,\tfrac23]. 
$$
For $N \ge 3$ the prefactor of $\alpha$ is negative and the lower bound \eqref{totcostbd} becomes minimal if and only if $\alpha$ becomes maximal, i.e. $\alpha=\tfrac23$, giving the total lower bound 
$$
  C[\gamma] \ge (4  + \tfrac{4}{N-1}) \cdot \tfrac13.
$$
The lower bound agrees with the cost \eqref{opticost} of the asserted  optimizer, completing the proof. 
\vspace*{2mm}

{\bf Proof of uniqueness (Prop.~\ref{P1}~c)).}  
Suppose $\gamma$ is an optimizer. Since equality must hold in \eqref{totcostbd} and $\alpha$ must be equal to $\tfrac23$, (i) all orbits starting at $\pm 1$ pass through zero and equality holds in \eqref{massineq} so that $\gamma(\mbox{\small orbits that vacate 0})=2/(3(N-1))$ and these orbits vacate $0$ for all $N-1$ timepoints $t_1,...,t_{N-1}$, (ii) each orbit that $\gamma$ gives mass to must achieve the optimal cost in \eqref{clb}, that is to say  $$
    \{\mbox{\small orbits that vacate 0}\}\cap {\rm supp}\,\gamma =\{(0,1,...,1,0),\, (0,\mi 1,...,\mi 1,0)\}
$$
and 
\begin{align*}
&\{\mbox{\small orbits that pass through 0}\}\cap{\rm supp}\,\gamma =\{\omega^{(i)}_+,\omega^{(i)}_-\}_{i=1}^{N-1} \\
&\, \mbox{ where }
\omega^{(i)}_{\pm} = \pm (1,...,1,0,-1,...,-1) \mbox{ with the value $0$ occuring at timepoint $t_i$}.
\end{align*}
So we have 
$$
  \gamma((0,1,...,1,0)) = \tfrac{1}{3(N-1)} + \delta, \;\; \gamma((0,\mi 1,...,\mi 1,0)) = \tfrac{1}{3(N-1)} - \delta \; \mbox{ for some }\delta\in[-\tfrac{1}{3(N-1)},\tfrac{1}{3(N-1)}]. 
$$
Mass conservation at $(x,t)=(\pm 1,t_1)$ implies 
$$
  \gamma((1,0,\mi 1,...,\mi 1))= \tfrac{1}{3(N-1)}+\delta, \;\;
  \gamma((\mi 1,0, 1,..., 1))= \tfrac{1}{3(N-1)}-\delta.
$$
Now using mass conservation at $(x,t)=(\pm 1,t_2)$ gives 
$$
  \gamma((1,1,0,\mi 1,...,\mi 1))= \tfrac{1}{3(N-1)}-\delta, \;\;
  \gamma((\mi 1, \mi 1, 0, 1,..., 1))= \tfrac{1}{3(N-1)}+\delta
$$
and by iteration we obtain 
\begin{equation} \label{mass_omega}
  \gamma(\omega^{(i)}_+) = \tfrac{1}{3(N-1)} + (-1)^{i-1}\delta, \;\;\;  \gamma(\omega^{(i)}_-) = \tfrac{1}{3(N-1)} - (-1)^{i-1}\delta  \;\;\;\; (i=1,...,N - 1).
\end{equation}
Now a difference between $N$ even and $N$ odd appears. When $N$ is even, the number of orbits $\omega_+^{(i)}$ ($i=1,...,N\mi 1$) is odd and eq.~\eqref{mass_omega} implies 
$$
    \gamma(\mbox{\small orbits that pass through 0 starting from 1}) = \sum_{i=1}^{N-1} \gamma(\omega_+^{(i)}) = \tfrac{N-1}{3(N-1)} + \delta,
$$
which together with $\gamma(\mbox{\small orbits that pass through 0 starting from 1})=\pi_0{}_\sharp \gamma(1)=\tfrac13$ implies $\delta=0$. This  shows that $\gamma_0$ is the unique optimizer. On the other hand, when $N$ is odd, eq.~\eqref{mass_omega} yields 
$$
    \gamma(\mbox{\small orbits that pass through 0 starting from 1}) = \sum_{i=1}^{N-1} \gamma(\omega_+^{(i)}) = \tfrac{N-1}{3(N-1)},
$$
so no restriction on $\delta$ arises from $\pi_0{}_\sharp \gamma(1)=\tfrac13$ and the set of optimizers is a one-parameter family, with different choices of $\delta$ giving different optimizers. For instance, $\delta=\tfrac{1}{3(N-1)}$ means that as compared to the plan $\gamma_0$ from Figure 2, mass initially at $1$ respectively $-1$ is moved in twice as large pieces through $0$, but at alternating times. This completes the proof of the proposition. 
\\[2mm]
Finally we note that when $N=3$, the new optimizer with $\delta=\tfrac{1}{3(N-1)}$ found in the proof of c$\mbox{\small )}$ is of Monge form, showing that the restriction $N\ge 4$ in Proposition \ref{P1} a$\mbox{\small )}$ is necessary. 

\section{Refining the spatial mesh}

Numerical simulations suggest that the mass-splitting phenomenon persists if the spatial mesh is refined, but the behavior of the solutions becomes somewhat complicated. Figure 3 below shows accurate numerical optimizers, obtained by solving the fully discrete generalized Euler equations (which are a linear program) with Matlab's inbuilt LP solver {\tt linprog}.\footnote{For a larger number of timesteps this approach quickly becomes infeasible since the number of unknowns increases exponentially with the number of timesteps. Tackling this curse of dimensionality requires much more sophisticated methods \cite{FP23, FP26}.} A rigorous explanation of mass-splitting in the spirit of Proposition \ref{P1} a) for general meshes would be desirable but lies beyond our scope. The delicacy of this question is illustrated by the occasional appearance of Monge solutions (see the top right panel). 

\begin{figure}
    \begin{center}
        \includegraphics[width=0.95\textwidth]{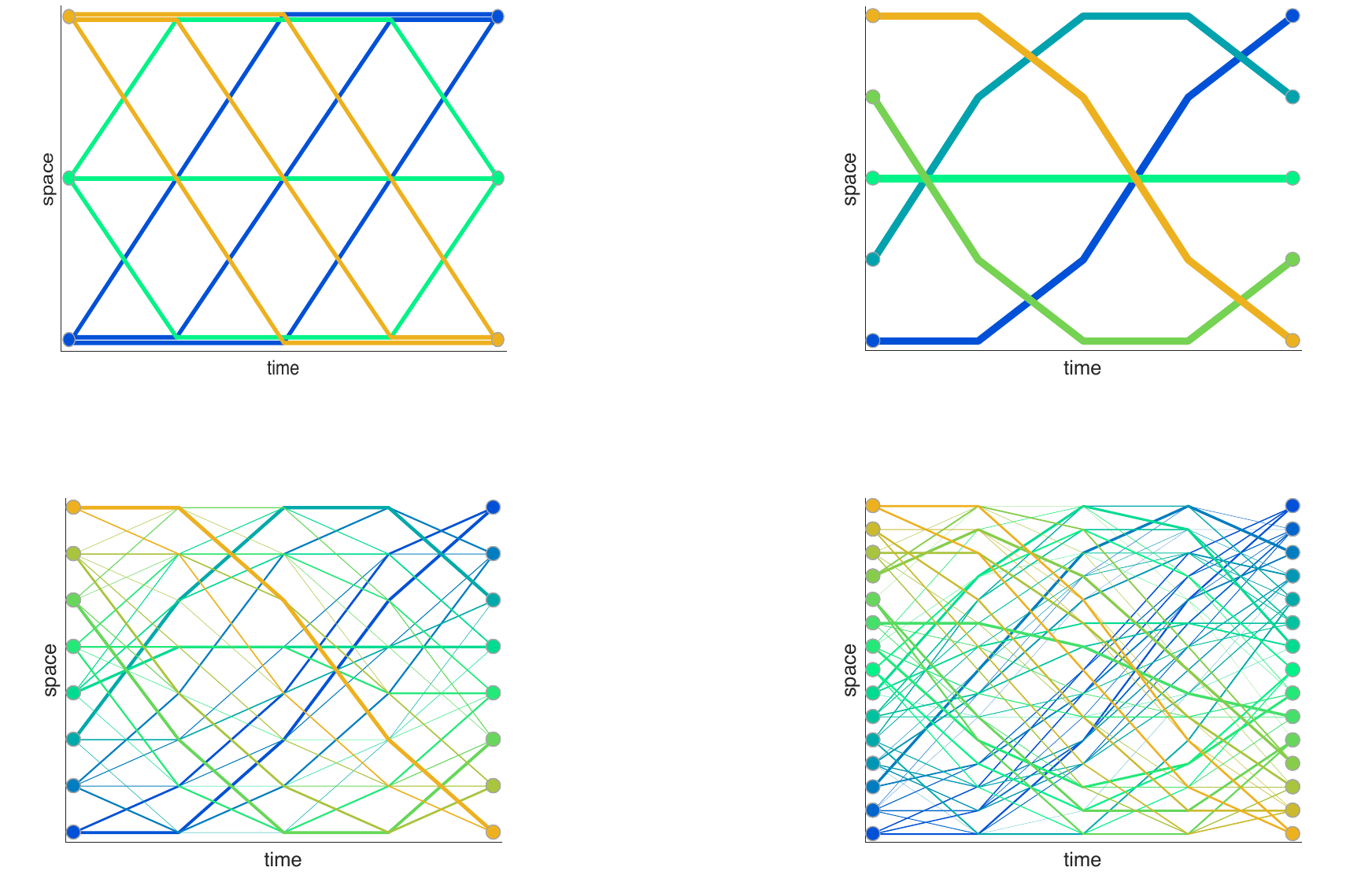}
    \end{center}
    \vspace*{-4mm}
    
    \caption{Numerical solutions of the generalized Euler equations on a space-time mesh. The thickness of the paths indicates the amount of mass transported.}
\end{figure}


\section{Example of mass-splitting on continuous space}
The time-discrete, spatially \textit{continuous}  generalized Euler equations \eqref{MMOT'1}--\eqref{MMOT'2} are a prototype example of multi-marginal optimal transport (MMOT), yet the question whether it exhibits mass-splitting has remained open. Here we answer this question positively.

\begin{theorem} \label{T1} Let $\Omega=[-1,1]$, $g_*(x)=-x$ (i.e. one turns the fluid upside down), $N=3$ (i.e. one considers three timesteps), $\{t_0,t_1,...,t_N\}=\{0,1,2,3\}$. The mass-splitting plan visualized in Figure 4, which is concentrated on the paths
$$
  \omega_{x,\rmv} \, : \, \{0,1,2,3\}\to[-1,1], \; \omega_{x,\rmv}(n) = x \cos \tfrac{\pi}3 n + \rmv \sin \tfrac{\pi}3 n 
$$
with probability distribution
\begin{align*} 
   f(x,\rmv) = \tfrac12 1_{[-1,1]}(x) \Bigl[\tfrac12  \delta_{\frac{3|x|-2}{\sqrt{3}}}(\rmv) + \tfrac12 \delta_{-\frac{3|x|-2}{\sqrt{3}}}(\rmv) \Bigr] \, dx \, d\rmv,
\end{align*}
is an optimizer of \eqref{MMOT'1}--\eqref{MMOT'2}. 
\end{theorem}

\begin{figure}[http!]
\begin{center}
    \includegraphics[width=0.65\textwidth]{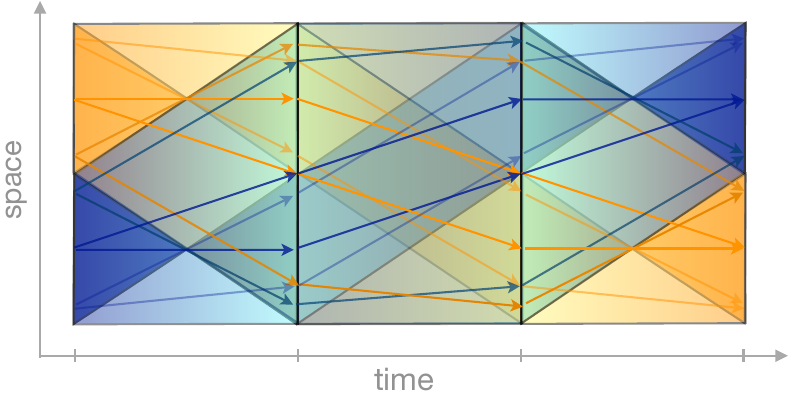}
\caption{Mass-splitting solution of the time-discrete generalized Euler equations \eqref{MMOT'1}--\eqref{MMOT'2}. Each fluid particle moves from its initial position to its end position via two 
paths.
}
\label{F:cts_example}
\end{center}
\end{figure}

The construction of this optimizer builds on the fully discrete example from Proposition \ref{P1} for $\Omega=\{-1,0,1\}$ and $N=3$ which is also concentrated
on two paths per initial position, 
$$
\begin{array}{rlrl}
   \omega\!\!\! &= (1,1,0,-1)  \mbox{ and } (1,0,-1,-1) \; \mbox{ for }x=1, \\
   \omega\!\!\! &= (0,1,1,0)  \;\;\,\mbox{ and }  (0,-1,-1,0) \;\mbox{ for }x=0, \\
   \omega\!\!\! &= (-1,0,1,1)  \mbox{ and }  (-1,-1,0,1) \;\mbox{ for }x=\!- 1.
\end{array}
$$
In fact these paths from Proposition \ref{P1} coincide \textit{exactly} with those used by the continuous optimizer from Theorem \ref{T1} for $x=1$, $0$, $-1$, as the reader can easily check using the explicit values of sine and cosine at multiples of $\tfrac{\pi}3$. What is more, the initial slopes of the two paths with initial position $x$ used by the continuous optimizer vary linearly with $x$ in $[-1,0]$ and $[0,1]$. Hence the continuous optimizer can be viewed as a certain piecewise linear interpolation of the discrete optimizer. 
It would be very interesting if such a construction  (first introduced by the author \cite{Fri19} in the context of a different MMOT problem) could be established in greater generality. 

The form of the paths in Theorem \ref{T1} comes from solving a time-discretized version of the Euler equations in Lagrangian coordinates, \eqref{concreteEL}, 
\begin{equation} \label{discreteEL}
    \omega(i\,\pl\, 1)-2\omega(i) + \omega(i\,\mi\, 1) \; =\; - \, p'\bigl(\omega(i)\bigr) \;\;\;\; (i=1,...,N-1), 
\end{equation}
with the explicit pressure 
\begin{equation} \label{pressure}
    p(x)=\tfrac12 x^2.
\end{equation}
The general solution to this equation is given precisely by the functions $\omega_{x,\rmv}$ in the theorem with $x\in\R$ and $\rmv\in\R$. For a derivation of \eqref{discreteEL}--\eqref{pressure} as a necessary condition on paths in the support of the optimal plan (and of its general solution) see the proof of the theorem. 

\begin{figure}[http!]
\begin{center}
    \includegraphics[width=0.65\textwidth]{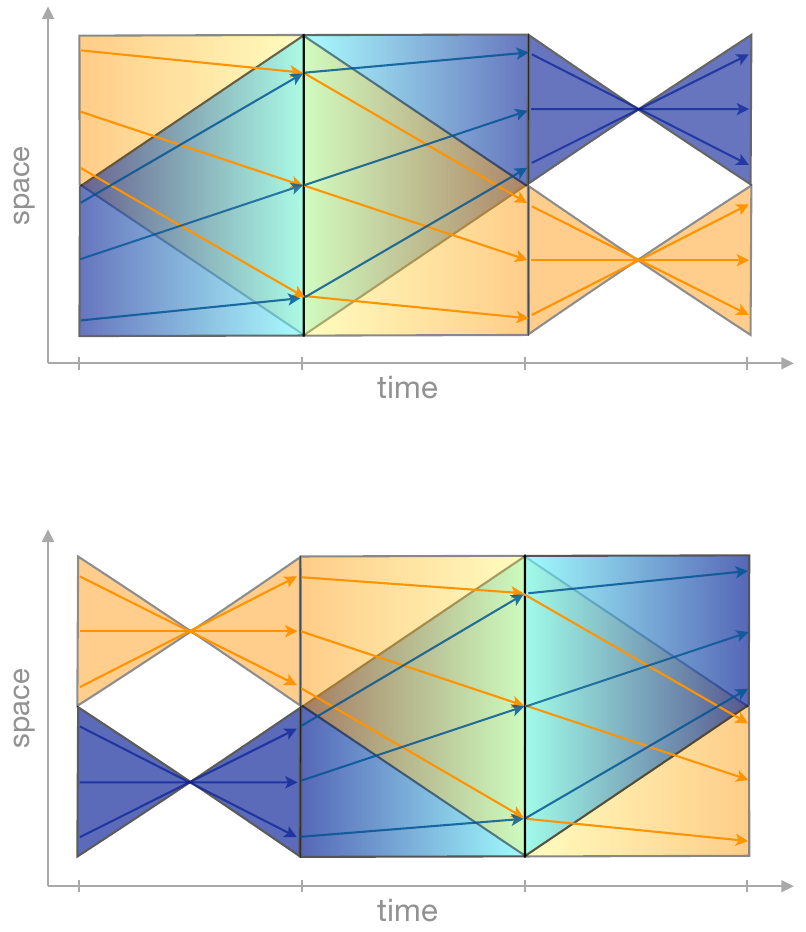}
\caption{The building blocks $\gamma_1$ (top) and $\gamma_2$ (bottom) of the optimal plan from Fig.~\ref{F:cts_example}. The top plan can be described as ``\textit{expand and mix, contract and de-mix, flip each block}''. The bottom plan corresponds to ``\textit{flip each block, expand and mix, contract and de-mix}''.
}
\label{F:cts_blocks}
\end{center}
\end{figure}

Our presentation of the example in terms of trigonometric functions has the advantage of making the connection with the Euler equations transparent (see \eqref{discreteEL}), but mass conservation appears a little mysterious. To understand the latter, it is helpful to decompose the optimal plan $\gamma_0$ into two simpler plans, $\gamma_0 = \tfrac12(\gamma_1 + \gamma_2)$ where $\gamma_1$ and $\gamma_2$ correspond to 
\begin{align*}
   & f_1(x,\rmv) = \tfrac12 1_{[-1,1]}(x) \left\{ \!\! \begin{array}{lr}\delta_{\frac{3|x|-2}{\sqrt{3}}}(\rmv) \!\!\! & \mbox{ if }x\ge 0 \\
   \delta_{-\frac{3|x|-2}{\sqrt{3}}}(\rmv) \!\!\! & \mbox{ if }x<0 \end{array} \!\right\} \, dx \, d\rmv, \\
   & f_2(x,\rmv) = \tfrac12 1_{[-1,1]}(x) \left\{ \!\! \begin{array}{lr}\delta_{-\frac{3|x|-2}{\sqrt{3}}}(\rmv) \!\!\! & \mbox{ if }x\ge 0 \\
   \delta_{\frac{3|x|-2}{\sqrt{3}}}(\rmv) \!\!\! & \mbox{ if }x<0 \end{array} \!\right\} \, dx \, d\rmv,
\end{align*}
and describe the underlying paths in a more elementary (trigonometric-function-free) manner. The plans $\gamma_1$ and $\gamma_2$ and the underlying paths $(\omega_0,...,\omega_3)$ are depicted in Figure \ref{F:cts_blocks}.
The paths correspond to certain piecewise linear measure-preserving maps $\omega_0\mapsto \omega_i$; the formulas are easy to read off from the figure and are given in the proof. 

As shown below, the plans $\gamma_1$ and $\gamma_2$ are themselves optimizers of \eqref{MMOT'1}--\eqref{MMOT'2}. These plans are only mass-splitting at the discrete time $t_1$ respectively $t_2$, but not at $t_0$. 
We do not know whether \textit{every} optimizer is mass-splitting at \textit{some} discrete time $t_i$. Our work only shows that \textit{some} optimizer is mass-splitting at \textit{every} $t_i$.  
\vspace*{2mm}

{\bf Proof of Theorem \ref{T1}.} \textit{Step 1.} The variational problem \eqref{MMOT'1}--\eqref{MMOT'2} under study is, explicitly, 
\vspace*{-3mm}
\begin{align} \label{MMOT'1ex}
&   \min_{\gamma\in\calP([-1,1]^4)}\, \int_{[-1,1]^4} \sum_{i=1}^3 |\omega_i-\omega_{i-1}|^2\, d\gamma(\omega_0,...,\omega_3) \; \mbox{ subject  to } 
 (\pi_0,\pi_3)_\sharp \gamma = (id,-id)_\sharp {\mathbf 1}_{[-1,1]} \\ 
 & \;\;\mbox{ and subject to the marginal conditions }\pi_i{}_\sharp \gamma = {\mathbf 1}_{[-1,1]}\mbox{ for } i=0,1,2,3, \label{MMOT'2ex}
\end{align}
where $\pi_i \, : \, \omega=(\omega_0,...,\omega_3)\mapsto \omega_i$ is the projection map onto the $i$-th coordinate. We will also make use of the reduced formulation \eqref{MMOT1}--\eqref{MMOT2} which in our case reads
\begin{align}
&   \min_{\gammabar\in\calP([-1,1]^{3})}\, \int_{[-1,1]^{3}} \Bigl( |\omega_0-\omega_1|^2 + |\omega_1-\omega_2|^2 + |\omega_2-(-\omega_0)|^2\Bigr)\, d\gammabar(\omega_0,\omega_1,\omega_2) \label{MMOT1ex} \\
& \; {\rm subject \; to} \; \pi_i{}_\sharp \gammabar = {\mathbf 1}_{[-1,1]} \, \mbox{ for }i=0,1,2. \label{MMOT2ex}
\end{align}
Recall from Lemma \ref{L:equiv} that a plan $\gamma$ is an optimizer of \eqref{MMOT'1ex}--\eqref{MMOT'2ex} if and only if its push-forward $(\pi_0,\pi_1,\pi_2)_\sharp \gamma$ is an optimizer of \eqref{MMOT1ex}--\eqref{MMOT2ex}. 
In the following we denote the projection of a path $\omega=(\omega_0,...,\omega_3)\in[-1,1]^4$ onto its first 3 components, $(\pi_0,\pi_1,\pi_2)(\omega)= (\omega_0,\omega_1,\omega_2)$, by $\omegabar$.

Next, we change the cost function $c(\omegabar)=|\omega_0-\omega_1|^2+|\omega_1-\omega_2|^2 + |\omega_2-(-\omega_0)|^2$ in \eqref{MMOT1ex} to 
\begin{equation} \label{modcost}
     \ctilde(\omegabar) = c(\omegabar) - (\omega_0^2+\omega_1^2+\omega_2^2).
\end{equation}
This does not change optimal plans, because integrating the extra terms with respect to any admissible plan for \eqref{MMOT1ex}--\eqref{MMOT2ex} only depends on the marginals, 
$$
\int \ctilde \, d\gammabar = \int c \, d\gammabar - \sum_{i=0}^2 \int \omega_i^2 \mathbf{1}_{[-1,1]}(\omega_i) d\omega_i.
$$

\textit{Step 2.} We now determine when the modified cost function $\ctilde$ is minimized pointwise. This function has the form $\ctilde(\omegabar) = \langle\omegabar,A\omegabar\rangle$ where $\langle \cdot , \cdot\rangle$ is the euclidean inner product on $\R^3$ and
$$
   A = \left( \!\!\! \begin{array}{rrr}
   1 & -1 & 1 \\ -1 & 1 & -1 \\ 1 & -1 & 1
   \end{array}\right).
$$
This is a rank-1 matrix, 
$$
    A = \left(\!\!\!\begin{array}{rrr} 1 \\ -1 \\ 1 \end{array}\right) \otimes \left(\!\!\!\begin{array}{rrr} 1 \\ -1 \\ 1 \end{array}\right)
$$
(where $a\otimes b$ denotes the matrix $M$ with components $M_{ij}=a_ib_j$). Thus 
$$
   \ctilde(\omegabar)\ge 0, \; \mbox{``=''} \Longleftrightarrow \Bigl\langle \omegabar, \left(\!\!\!\begin{array}{rrr} 1 \\ -1 \\ 1 \end{array}\right) \Bigr\rangle = 0.
$$
Moreover if $\omega=(\omegabar,\omega_3)$ satisfies the condition on the right, that is to say 
$$
      \omega_0-\omega_1+\omega_2 =0,
$$
then we have $\omega_1-\omega_2+\omega_3=0$ if and only if $\omega_3=-\omega_0$. Hence the following equivalence holds for any path $\omega=(\omegabar,\omega_3)\in\R^4$: 
\begin{align*}
    \ctilde(\omegabar) = \min \ctilde \; \mbox{ and } \;  \omega_3=\omega_0 \;\;  &\Longleftrightarrow 
    \omega_{i-1}-\omega_i+\omega_{i+1}=0 \; \mbox{ for }i=1,2 \\
    &\Longleftrightarrow 
    \omega \mbox{ satisfies \eqref{discreteEL}--\eqref{pressure}.}
\end{align*}
Here the last equivalence follows trivially by adding $(-\omega_i)$ to both sides of the equation and noting that $N=3$. 
\vspace*{2mm}

\textit{Step 3.} The solutions to \eqref{discreteEL}--\eqref{pressure} are given by the following lemma which is well-known from the theory of the discrete Laplacian in one dimension. 
\begin{lemma} \label{L:trig} For any $N\ge 2$, the general solution $\omega = (\omega_0,...,\omega_N)$ to \eqref{discreteEL} is 
\begin{equation} \label{solformula}
  \omega_n = x \cos \tfrac{n\pi}3  + \rmv \sin \tfrac{n\pi}3 \;\; (n=0,...,N)
\end{equation}
\vspace*{-7mm}

\noindent
with $x$, $\rmv\in\R$. 
\end{lemma}
\noindent
{\bf Proof} (included for convenience of the reader). Write the second order difference equation \eqref{discreteEL} as a recursion 
$$
   \begin{pmatrix} \omega_n \\ \omega_{n+1} \end{pmatrix} = \underbrace{\begin{pmatrix} 0 & 1 \\ -1 & 1 \end{pmatrix}}_{=:B} \begin{pmatrix} \omega_{n-1} \\ \omega_n\end{pmatrix}, \;\;\; n=1,2,...
$$
so that
\begin{equation*} 
   \begin{pmatrix} \omega_n \\ \omega_{n+1} \end{pmatrix} = B^n \begin{pmatrix} \omega_{0} \\ \omega_1\end{pmatrix}.
\end{equation*}
The matrix $B$ has eigenvalues $\lambda_\pm = \tfrac{1\pm i\sqrt{3}}{2}=e^{\pm i\pi/3}$ and corresponding eigenvectors 
$$
   v_\pm = \begin{pmatrix} 1 \\ \lambda_\pm \end{pmatrix}. 
$$
Decomposing the initial vector  $(\omega_0,\omega_1)$ of the recursion as $(\omega_0,\omega_1)=\alpha v_+ + \beta v_-$ for some coefficients $\alpha$, $\beta\in\C$ gives
$$
   \begin{pmatrix} \omega_n \\ \omega_{n+1} \end{pmatrix} = \alpha \lambda_+^n v_+ + \beta \lambda_-^n v_-
$$
and thus
$$
   \omega_n = \alpha e^{i\tfrac{n\pi}3 } + \beta e^{-i\tfrac{n\pi}3 }. 
$$
Passing to real trigonometric functions gives the representation in the lemma.
\vspace*{2mm}

For physical interpretation, let us express the coefficients $x$ and $\rmv$ in terms of $\omega_0$ and $\omega_1$. Evaluating \eqref{solformula} at $n=0$ gives 
\begin{equation} \label{xeq}
   x = \omega_0,
\end{equation}
so $x$ is the initial position. 
Extending the recursion backwards via $\omega_{-1}-\omega_0+\omega_1=0$, i.e. setting $\omega_{-1}:=\omega_0-\omega_1$, and using $\omega_{\pm 1}= x \cos \tfrac{\pi}3 \pm \rmv \sin \tfrac{\pi}3$ gives
$$
   \omega_1-\omega_{-1} = 2 \rmv \sin\tfrac{\pi}3 = 2 \rmv \tfrac{\sqrt{3}}{2}
$$
and so 
\begin{equation} \label{veq}
   \rmv = \frac{\omega_1-\omega_{-1}}{\sqrt{3}}.
\end{equation}
Thus $\rmv$ is a discrete initial velocity, defined via a central difference. 

Now let $\gamma_0$ be any plan supported on paths $\omega=(\omegabar,\omega_3)$ satisfying \eqref{discreteEL}--\eqref{pressure} (as is the case for the plan given in the theorem, as well as for the auxiliary plans $\gamma_1$ and $\gamma_2$ introduced below the theorem). By Step 2 and Lemma \ref{L:trig}, all paths in the support satisfy $\ctilde(\omegabar)=\min \ctilde$ and so the push-forward $(\pi_0,\pi_1,\pi_2)_\sharp\gamma_0$ is a minimizer of $\int \ctilde d\gammabar$ over all of $\calP([-1,1]^3)$. Thus by Step 1, provided $\gamma_0$ satisfies the marginal conditions \eqref{MMOT'2ex} it is a minimizer of \eqref{MMOT'1ex}--\eqref{MMOT'2ex}.   
\vspace*{2mm}

\textit{Step 4.} Finally we show that the plan $\gamma_0$ given in the theorem satisfies the marginal conditions \eqref{MMOT'2ex}. By the decomposition $\gamma_0=\tfrac12(\gamma_1 + \gamma_2)$ introduced below the theorem, it suffices to check the marginal conditions for $\gamma_1$ and $\gamma_2$ (thereby showing that these plans are also minimizers of \eqref{MMOT'1ex}--\eqref{MMOT'2ex}). We use that 
$$
 \cos \tfrac{\pi}3 = \tfrac12, \;\; \cos \tfrac{2\pi}{3}=-\tfrac12, \;\; \sin \tfrac{\pi}3 = \sin \tfrac{2\pi}{3}=\tfrac{\sqrt{3}}{2}. 
$$
Hence paths in the support of $\gamma_1$ satisfy, using \eqref{xeq}, 
\begin{align*}
  \omega_1 & = \omega_0  \cos \tfrac{\pi}3 
  + \left\{ \!\! \begin{array}{lr} \frac{3|\omega_0|-2}{\sqrt{3}} \!\!\! & \mbox{ if }\omega_0\ge 0 \\
    \frac{-3|\omega_0|+2}{\sqrt{3}} \!\!\! & \mbox{ if }\omega_0<0 \end{array} \!\right\}\sin \tfrac{\pi}3
  = \frac{\omega_0}2 + \begin{cases} \tfrac{3|\omega_0|-2}{2} & \mbox{if }\omega_0\ge 0 \\
  \tfrac{-3|\omega_0|+2}{2} & \mbox{if }\omega_0 < 0 \end{cases} \\
  & =  \begin{cases} 2\omega_0 - 1 \!\!\! & \mbox{ if }\omega_0\ge 0 \\
   2\omega_0 + 1 \!\!\! & \mbox{ if }\omega_0<0\end{cases} \; =: T_1(\omega_0).
\end{align*}
This is the ``\textit{expand and mix}'' map visualized in the top left part of Figure \ref{F:cts_blocks}. 
Further, 
$$
  \omega_2 \underset{\eqref{discreteEL},\,\eqref{pressure}}{=} \omega_1-\omega_0 = 
   \begin{cases} \omega_0 - 1 \!\!\! & \mbox{ if }\omega_0\ge 0 \\
   \omega_0 + 1 \!\!\! & \mbox{ if }\omega_0<0\end{cases} \; =: T_2(\omega_0).
$$
Thus $\gamma_1=(id,T_1,T_2,-id)_\sharp \mathbf{1}_{[-1,1]}$, and the marginal conditions \eqref{MMOT'2ex} follow since $T_1$, $T_2$ are obviously measure-preserving maps on $[-1,1]$. 
Analogously, the paths in the support of $\gamma_2$ satisfy 
   \begin{align*}
  \omega_1 & = \frac{\omega_0}2 + \begin{cases} \tfrac{-3|\omega_0|+2}{2} & \mbox{if }\omega_0\ge 0 \\
  \tfrac{3|\omega_0|-2}{2} & \mbox{if }\omega_0 < 0 \end{cases}  \; = \;  \begin{cases} -\omega_0 + 1 \!\!\! & \mbox{ if }\omega_0\ge 0 \\
   -\omega_0 - 1 \!\!\! & \mbox{ if }\omega_0<0\end{cases} \; =: S_1(\omega_0)
\end{align*}
(this is the ``\textit{flip each block}'' map visualized in the bottom left part of Figure \ref{F:cts_blocks}) and
$$
  \omega_2 \; = \; \omega_1-\omega_0 \, = 
   \begin{cases} -2\omega_0 + 1 \!\!\! & \mbox{ if }\omega_0\ge 0 \\
   -2\omega_0 - 1 \!\!\! & \mbox{ if }\omega_0<0\end{cases} \; =: S_2(\omega_0).
$$
Thus $\gamma_2=(id,S_1,S_2,-id)_\sharp {\mathbf 1}_{[-1,1]}$, and the marginal conditions \eqref{MMOT'2ex} follow also for $\gamma_2$ since $S_1$, $S_2$ are also measure-preserving. The proof of the theorem is complete.

\section{What is physically more correct, Euler (no mass splitting) or Brenier (mass splitting)?}
\label{sec:modeling}

We close with some remarks from a modeling point of view. 

Continuum mechanical models of fluids, like the Euler equations, are macroscopic descriptions of the collective motion of the underlying microscopic particles. In the case of water, these are H$_2$O molecules. Neglecting quantum effects 
(such as
the fact that these molecules occasionally split into their constituents H, OH, H$^+$ and OH$^-$ and recombine) leads to molecular dynamics models. These describe the system by the positions $X_j$ and velocities $V_j$ of the particles, and typically take the form of an evolution equation $m_i\ddot{X}_i(t) = f_i(\{X_j\},\{V_j\})$ (where $m_j$ is the mass of the $j$-th particle), or equivalently 
\begin{equation} \label{MD}
   \begin{array}{l}\dot{X}_i = V_i, \\[1mm]
   m_i\dot{V_i} = f_i(\{X_j\},\{V_j\}).
   \end{array}
\end{equation}

Mesoscopic models describe the state of the system at time $t$ by a phase space density $f(x,\rmv,t)$ of particles at position $x$ with velocity $\rmv$ at time $t$.  
For very dilute systems (rarefied gases), an accurate  mesoscopic model for the time evolution is given by the Boltzmann equation
\begin{align} \label{Boltzmann}
   \partial_t f + ({\rm v}\cdot \nabla_x)\, f = Q(f,f)
\end{align}
where $Q$ is a quadratic interaction kernel. 

Continuum mechanical models coarse-grain the phase space density even further. Typically, one describes the system by a position density $\rho(x,t)$ of particles at position $x$ at time $t$ and a single velocity $u(x,t)$ at position $x$ at time $t$ (representing an average particle velocity), and one puts forth a system of evolution equations for $\rho$ and $u$. In the case of the Euler equations for an incompressible fluid with constant density $\rho$, these can be written in the form 
\begin{equation} \label{contmech}
   \begin{array}{l} 
   \partial_t \rho + {\rm div}(\rho u)=0, \\[1mm]
   \partial(\rho u) + {\rm div}(\rho u \otimes u) = -\nabla p.
   \end{array}
\end{equation}

It is instructive\footnote{nonwithstanding the fact that the Boltzmann equation is valid for a low-density gas, whereas the Euler equations model a high-density  system like liquid water} to compare Euler's system with the (non-closed) evolution system for density and average velocity that can be extracted from the Boltzmann equation \eqref{Boltzmann}. In terms of Boltzmann's phase space density $f$, density and average velocity are 
$$
   \rho(x,t) = \int f(x,{\rm v},t)\, d{\rm v}, \;\;\;\; u(x,t) = \frac{\int {\rm v}\, f \, d{\rm v}}{\int f \, d{\rm v}}.
$$
Multiplying the equation by $1$ and integrating over $\rmv$ gives the density evolution:
\begin{align*}
 & \partial_t \int f \, d{\rm v} + \nabla_x\cdot \int {\rm v} \, f \, d{\rm v} = 0 
\end{align*}
or equivalently 
\begin{align*}  
 \partial_t \rho + {\rm div}(\rho u) = 0.
\end{align*}
This agrees with the first equation in \eqref{contmech}. Multiplying the Boltzmann equation by $\rmv$ and integrating over $\rmv$ gives the velocity evolution: 
\begin{align*}
    \partial_t(\rho \, u) = {\rm terms \; depending \; on \; higher \; moments \; of} \; f \; \mbox{\rm  with respect to  v}.
\end{align*}
This does not agree with the Euler equation (second equation in \eqref{contmech}), which ignores higher moments (velocity fluctuations).

We find it noteworthy that Brenier's relaxation of Euler brings certain velocity fluctuations back in, as observed in the pioneering paper \cite{Bre89} and illustrated further by the new examples in this paper. 
The velocity distributions seen in generalized Euler are very different from the Maxwellian distributions emerging at long time from the Boltzmann equation. It would be interesting to attempt to derive such distributions, or at least a selection thereof, from  concrete microscopic or mesoscopic models. 



\vspace*{4mm}

\noindent
{\bf Acknowledgements.} The author thanks Yann 
Brenier, Maximilian Penka and Luca Nenna for helpful discussions on the generalized Euler equations. 

\begin{small}
\begin{spacing}{0.9}



\end{spacing}

\end{small}

\end{document}